\documentclass{article}
\author{T. S. Trudgian\footnote{Supported by ARC Grant DE120100173.}\\
Mathematical Sciences Institute\\ The Australian National University,
 ACT 0200, Australia\\ timothy.trudgian@anu.edu.au
}

\title{A short extension of two of Spira's results}
\usepackage{url}
\usepackage{amsthm}
\usepackage{amsmath}
\usepackage{comment}
\usepackage{amssymb}
\usepackage{booktabs}
\newtheorem{thm}{Theorem}

\newtheorem{cor}{Corollary}
\begin{document}

\maketitle

\begin{abstract}
\noindent
Two inequalities concerning the symmetry of the zeta-function and the Ramanujan $\tau$-function are improved through the use of some elementary considerations.
\end{abstract}
\section{Introduction}
The functional equation for the Riemann zeta-function $\zeta(s)$ is 
$$\zeta(1-s) = g(s) \zeta(s), \quad\quad g(s) = 2^{1-s} \pi^{-s} \cos \tfrac{1}{2} s\pi \Gamma(s),$$ 
see \cite[(2.1.8)]{Titchmarsh}.
Spira \cite{Spira65} proved that 
\begin{equation}\label{spire}
|\zeta(1-s)| > |\zeta(s)|, \quad \tfrac{1}{2} < \sigma <1,\quad t\geq 10, \quad \textrm{when}\; \zeta(s)\neq 0.
\end{equation}
Dixon and Schoenfeld \cite{DixonS} gave a simpler and sharper proof of (\ref{spire}) for $|t|\geq 6.8$ and for all $\sigma >\frac{1}{2}$. Saidak and Zvengrowski \cite{Saidak} proved (\ref{spire}) for $|t|\geq 2\pi +1$, and, in fact, their proof is valid for $|t|\geq 7$. 
Recently, Nazardonyavi and Yakubovich\footnote{The authors considered the equivalent problem of $|\zeta(1-s)| < |\zeta(s)|$ for $0< \sigma < \frac{1}{2}$.} \cite{NazYak} gave an alternative proof of (\ref{spire}) in the range $|t|\geq 12$. They remark that this result may be extended to $|t|\geq 6.5$ by a computer simulation.
Spira [\textit{op.\ cit.}] notes that (\ref{spire}) `fails for $t$ around $2\pi$'.
Indeed, for $t^{*} = 6.2898$ one may compute
\begin{equation}\label{false}
\frac{|\zeta(0.48 - it^{*})|}{|\zeta(0.52 + it^{*})|} - 1 < -8\times 10^{-8}.
\end{equation}

The purpose of this short article is to examine the proof given by Dixon and Schoenfeld and to prove
\begin{thm}\label{The}
$|\zeta(1-s)| > |\zeta(s)|$ except at the zeroes of $\zeta(s)$, where $|t|\geq 6.29073$ and $\sigma >\frac{1}{2}$.
\end{thm}
By the functional equation, $\zeta(1-s)$ and $\zeta(s)$ have the same zeroes when $0<\sigma <1$. This means that Theorem \ref{The}  gives rise to the following Corollary, which improves on Proposition 1 in \cite{NazYak}.
\begin{cor}\label{only}
A necessary and sufficient condition for the Riemann hypothesis is
$$|\zeta(1-s)| > |\zeta(s)|, \quad  \sigma > \tfrac{1}{2},\quad |t| \geq 6.29073.$$
\end{cor}
In light of (\ref{false}) Theorem \ref{The} and Corollary \ref{only} are close to best possible. The purpose of this article is to show that the range of $t$ in (\ref{spire}) can be extended relatively easily. In  \cite{NazYak} the range is increased at the cost of significant computation. By contrast, almost no computation is required to establish Theorem \ref{The} and Corollary \ref{only}.

Similarly, for $F(s) = \sum_{n=1}^{\infty} \tau(n)/n^{s}$, where $\tau(n)$ is the Ramanujan $\tau$-function, Spira \cite{Spira73} proved that 
\begin{equation}\label{Fdef}
|F(12-s)|> |F(s)|, \quad 6<\sigma < \frac{13}{2}, \quad t\geq 4.35,
\end{equation}
except at the zeroes of $F(s)$. This improved on a result of Berndt \cite{Berndt70} who proved (\ref{Fdef}) for $t\geq 6.8$. At no extra charge, the proof of Theorem \ref{The} gives
\begin{thm}\label{t2}
$|F(12-s)|> |F(s)|$ except at the zeroes of $F(s)$, where $t\geq 3.8085$ and $6<\sigma < \frac{13}{2}$.
\end{thm}

\section{Proof of Theorems \ref{The} and \ref{t2}}
Dixon and Schoenfeld consider the function $h(s) = \log|g(s)/g(\frac{1}{2} + it)|$. Since $|g(\frac{1}{2} + it)| =1$, one may prove Theorem \ref{The} by showing that $h(s)>0$ for $\sigma>\frac{1}{2}$ and $t\geq 6.29073$.

Starting at \cite[(1)]{DixonS} we have
\begin{equation*}
\frac{h(s)}{\sigma - \frac{1}{2}} > \left\{ \frac{\partial}{\partial\sigma} \log| \Gamma(\sigma + it)|\right\}_{\sigma = \sigma_{1}} - 2\pi e^{-\pi t} - \log 2\pi,
\end{equation*}
for some number $\sigma_{1}\in(1/2, \sigma)$. Using Stirling's formula we arrive at \cite[(3)]{DixonS} which is
\begin{equation}\label{clarke}
\frac{\partial}{\partial\sigma} \log|\Gamma(\sigma + it)| = \Re\left\{ \log s - \frac{1}{2s} - \frac{1}{12 s^{2}} + 6 \int_{0}^{\infty} \frac{P_{3}(\{x\})}{(s+ x)^{4}}\, dx\right\},
\end{equation}
where $P_{3}(x) = x(2x^{2} - 3x + 1)/12$ and $\{x\}$ denotes the fractional part of $x$. Some simple calculus gives $\max_{x\in[0,1]} P_{3}(x) = \sqrt{3}/216$. 

Rather than bound each term in (\ref{clarke}) by its modulus, as in \cite{DixonS}, we consider each real part separately. For $s = \sigma + it$ the right-hand side of (\ref{clarke}) is bounded below by
\begin{equation}\label{his}
\frac{1}{2} \log(\sigma^{2} + t^{2}) - \frac{\sigma}{2(\sigma^{2} + t^{2})} -\frac{(\sigma^{2} - t^{2})}{12(\sigma^{2} + t^{2})^{2}} - \frac{\sqrt{3}}{36}\int_{0}^{\infty} \frac{dx}{\{(\sigma + x)^{2} + t^{2}\}^{2}}.
\end{equation}
The integral in (\ref{his}), denoted by $I$, is clearly decreasing in $\sigma$, whence we conclude
$$ I \leq \int_{0}^{\infty} \frac{dx}{\{(\frac{1}{2} + x)^{2} + t^{2}\}^{2}} = \frac{\tan^{-1} 2t - \frac{2t}{4t^{2} + 1}}{2t^{3}}.$$
Denote the first three terms in (\ref{his}) by $J(\sigma, t)$. It is easy to show that $$ \frac{\partial J}{\partial \sigma} = \frac{\sigma^{3}(1+ 3\sigma + 6\sigma^{2}) + 3t^{2}(\sigma - \frac{1}{2})\{ 2t^{2} + 4\sigma(\sigma + \frac{1}{2})\}}{6(\sigma^{2} + t^{2})^{3}},$$
which is clearly positive for $\sigma \geq \frac{1}{2}$, whence 
$$ \frac{h(s)}{\sigma - \frac{1}{2}} > J(\tfrac{1}{2}, t) - \frac{\sqrt{3}(\tan^{-1} 2t - \frac{2t}{4t^{2} + 1})}{72t^{3}}             - 2\pi e^{-\pi t} - \log 2\pi = G(t),$$
say.
It is straightforward to check that $G(t)$ is increasing and that $G(6.29072)<0 < G(6.29073)$, which proves Theorem \ref{The}.

To prove Theorem \ref{t2} we note that, by Spira \cite[p.\  384]{Spira73}, it is sufficient to show that
$$ \left\{\frac{\partial}{\partial \sigma}\log|\Gamma(\sigma + it)|\right\}_{\sigma = \sigma_{1}} - \log 2\pi>0,$$
where $\sigma_{1}\in[\frac{11}{2}, \frac{13}{2}]$. But for a small alteration in bounding the integral $I$ in (\ref{his}), the calculation proceeds as before. It is sufficient to show that
$$H(t) = J(\tfrac{11}{2}, t) - \frac{\sqrt{3}(\tan^{-1} \frac{2t}{11} - \frac{22t}{121 + 4t^{2}})}{72t^{3}} >0.$$
A computational check shows that $H(t)$ is increasing, and that $H(3.8024)<0<H(3.8085)$, which establishes Theorem \ref{t2}.

\subsection*{Acknowledgements}
I am grateful to Kevin Broughan and to Peter Zvengrowski for some fruitful discussions about this problem, and to the referee for a careful reading of the manuscript.

\bibliographystyle{plain}
\bibliography{themastercanada}

\begin{thebibliography}{1}

\bibitem{Berndt70}
B.~C. Berndt.
\newblock On the zeros of a class of {D}irichlet series. {I}.
\newblock {\em Illinois J. Math.}, 14:244--258, 1970.

\bibitem{DixonS}
R.~D. Dixon and L.~Schoenfeld.
\newblock The size of the {R}iemann zeta-function at places symmetric with
  respect to the point $\frac{1}{2}$.
\newblock {\em Duke Math. J.}, 33:291--292, 1966.

\bibitem{NazYak}
S.~Nazardonyavi and S.~Yakubovich.
\newblock Another proof of {S}pira's inequality and its application to the
  {R}iemann hypothesis.
\newblock {\em J. Math. Inequal.}, 7:167--174, 2013.

\bibitem{Saidak}
F.~Saidak and P.~Zvengrowski.
\newblock On the modulus of the {R}iemann zeta function in the critical strip.
\newblock {\em Math. Slovaca}, 53(2):145--172, 2003.

\bibitem{Spira65}
R.~Spira.
\newblock An inequality for the {R}iemann zeta function.
\newblock {\em Duke Math. J.}, 32:247--250, 1965.

\bibitem{Spira73}
R.~Spira.
\newblock Calculation of the {R}amanujan $\tau$-{D}irichlet series.
\newblock {\em Math. Comp.}, 27(122):379--385, 1973.

\bibitem{Titchmarsh}
E.~C. Titchmarsh.
\newblock {\em The Theory of the Riemann zeta-function}.
\newblock Oxford Science Publications. Oxford University Press, Oxford, 2nd
  edition, 1986.

\end{thebibliography}

\end{document}